\pdfoutput=1
\newif\ifarxiv
\arxivtrue   %
\ifarxiv
\documentclass[preprint,12pt]{elsarticle}
\else
\documentclass[review,12pt]{elsarticle}
\fi

\usepackage{graphicx}

\usepackage{tabularx}
\usepackage{amssymb}
\usepackage{amsmath}
\usepackage{float}
\usepackage{booktabs}
\usepackage{tabularx} %
\usepackage{array}    %

\usepackage{listings}
\usepackage{xcolor}

\lstdefinestyle{mystyle}{
    backgroundcolor=\color{black!5},
    commentstyle=\color{green!40!black},
    keywordstyle=\color{blue},
    numberstyle=\tiny\color{black!50},
    stringstyle=\color{purple},
    basicstyle=\ttfamily\footnotesize,
    breakatwhitespace=false,         
    breaklines=true,                 
    captionpos=b,                    
    keepspaces=true,                 
    numbers=left,                    
    numbersep=5pt,                  
    showspaces=false,                
    showstringspaces=false,
    showtabs=false,                  
    tabsize=2,
    language=Python
}
\usepackage{tikz}
\usepackage{tikz-3dplot}
\usepackage{microtype}
\usepackage{geometry}
\usepackage{orcidlink}

\newcounter{bla}

\ifarxiv
\myfooter[L]{\parbox[b]{\textwidth}{\raggedright Accepted manuscript. The version of record is published in Computer Physics Communications (2026) 110455, \textup{\url{https://doi.org/10.1016/j.cpc.2026.110455}}.}}
\fi
\begin{document}

\begin{frontmatter}

\title{PySymmetry: A Sage/Python Framework for the Symmetry Reduction of Linear G-Equivariant Systems}

\author[1]{Leon D. da Silva\orcidlink{0000-0002-2823-3418}}
\ead{leon.silva@ufrpe.br}

\author[1]{Marcelo P. Santos\orcidlink{0000-0002-9023-0728}}
\ead{marcelo.pedrosantos@ufrpe.br}

\address[1]{Federal Rural University of Pernambuco, Department of Mathematics, Dom Manoel de Medeiros Street, S/N, Dois Irmãos, Recife, PE, 52171-900, Brazil}

\begin{abstract}
Linear systems with finite-group symmetry, from geometric invariances to repetitive structural patterns, are ubiquitous in science and engineering, and standard computational software rarely exploits that structure. This paper introduces PySymmetry, an open-source Sage/Python framework that applies the constructive algorithms of classical representation theory to the symmetry reduction of $G$-equivariant linear systems. Projection operators built from the characters of the irreducible representations of a finite group generate a symmetry-adapted basis, which carries an equivariant operator into block-diagonal form and decouples one large problem into several small independent ones. A single character-theoretic construction serves two paths inside one open environment: an exact symbolic path over an exact ring, and a numerical path interoperable with \texttt{NumPy}, \texttt{SciPy}, and sparse matrix formats. Three case studies exercise both. The GF method of molecular vibrations illustrates the workflow. A non-Hermitian Schrödinger eigenvalue benchmark shows the reduced path to be nearly 17 times faster than a general-purpose dense eigensolver on the largest grid tested. A symbolic investigation produced the first complete analytical classification of the nested cube central configuration in celestial mechanics. Source code, documentation, and the notebooks that reproduce the results are publicly available.\\

\noindent \textbf{PROGRAM SUMMARY}

\begin{small}
	\noindent
	{\em Program Title: }\texttt{PySymmetry} \\
	{\em CPC Library link to program files:} \url{https://doi.org/10.17632/xm3skv4v68.1} \\
	{\em Developer's repository link:} \url{https://github.com/ldsufrpe/pysymmetry} \\
	{\em Licensing provisions:} MIT \\
	{\em Programming language:} SageMath 10.4, Python 3.11\\
	{\em Supplementary material:} The GitHub repository includes Jupyter notebooks that reproduce the figures and tables generated by the framework; the Software Availability section maps each one to the notebook that produces it. \\
	{\em Nature of problem:}
	Linear systems governed by operators with discrete symmetries ($G$-equivariant systems) are common in physics and engineering, from quantum mechanics to structural analysis. Solving these systems, such as large-scale eigenvalue problems ($Mx=\lambda x$) or linear equations ($Mx=b$), can be computationally prohibitive. In numerical contexts, the cost of standard dense eigensolvers scales poorly with system size. In symbolic contexts, the problem is often intractable due to intermediate expression swell, where symbolic calculations exhaust memory and processing resources. This framework addresses the need for a tool that can systematically exploit these inherent symmetries to reduce computational complexity in both numerical and symbolic domains.
	\\
	{\em Solution method:}
	The \texttt{PySymmetry} framework implements the constructive algorithms of group representation theory to simplify $G$-equivariant linear operators. Given a finite group $G$ and a representation, the method first uses characters and the conjugacy classes of the group to determine the structure of the representation's isotypic decomposition. It then algorithmically constructs projection operators to generate a symmetry-adapted basis for each isotypic component. By applying a similarity transformation with this basis, the original equivariant operator is transformed into a block-diagonal form. This procedure effectively decouples a single, large-scale problem into a set of smaller, independent subproblems that can be solved more efficiently, either numerically or symbolically.
	\\
	{\em Additional comments including restrictions and unusual features:}
	A key feature of \texttt{PySymmetry} is its hybrid symbolic-numeric architecture, built on the SageMath environment to provide seamless interoperability with the scientific Python ecosystem (NumPy, SciPy). It is designed with an object-oriented API and supports sparse matrix formats for efficiency. A significant advantage is the reusability of the computed symmetry-adapted basis: for parametric studies where the operator changes but the symmetry remains, the expensive preprocessing step is performed only once and its cost is amortized over numerous subsequent rapid calculations. The resulting block-diagonal structure is also inherently suited for parallel computation (``embarrassingly parallel''), further enhancing its performance potential on modern hardware.

\end{small}

\end{abstract}


\end{frontmatter}
\section{Introduction}

The application of group theory to simplify problems with symmetry is a practice that has been established for decades in physics and chemistry, and in recent decades, it has also been systematically exploited in engineering \cite{Zingoni2009} and in both pure and applied mathematics \cite{stiefel2012group, serre1977linear}. In this context, linear systems that exhibit symmetry under the action of finite groups arise in a wide range of scientific and engineering contexts \cite{Dresselhaus2007}, including quantum mechanics~\cite{castelazo2021quantum, kahnert2005irreducible}, quantum computing~\cite{castro2025symmetric}, solid-state physics~\cite{hergert2018group}, structural analysis~\cite{kaveh2013optimal}, celestial mechanics~\cite{santos2022decomposition, leandro2019structure, leandro2017factorization, santos2024symmetric, xia2022applying}, and boundary-value problems in computational mechanics~\cite{bossavit1986boundary,allgower1992exploiting,georg1994exploiting}. Classical problems such as solving linear systems \(Mx = b\) or eigenvalue problems \(Mx = \lambda x\), where the operator \(M\) is \(G\)-equivariant with respect to a finite group \(G\), possess intrinsic symmetries that can be exploited for computational gain~\cite{georg1994exploiting}. When properly exploited, these symmetries enable the transformation of the system into block-diagonal form, decoupling it into smaller, independent subsystems. As demonstrated in several works on the discretization of partial differential equations and boundary integral equations~\cite{allgower1992exploiting}, such symmetry-aware formulations yield meaningful decompositions that reduce computational complexity, improve memory efficiency, and enhance the accuracy of approximate solutions.

When such $G$-equivariant systems require exact symbolic solutions, however, the challenge shifts from one of mere efficiency to one of fundamental tractability. Here, the primary obstacle is the phenomenon of intermediate expression swell, where the size of intermediate expressions in an exact calculation grows exponentially, often rendering problems computationally intractable, not merely slow \cite{mansfield1997applications,vonZurGathen2013modern}. This barrier is particularly acute in symbolic linear algebra. Operations like computing determinants, inverses, or null spaces of the large matrices that model these symmetric structures can exhaust all available memory and processing resources long before a final result is obtained. Paradoxically, the very symmetry that defines these complex systems also offers the key to their simplification. Exploiting this inherent structure to simplify the matrices before performing costly symbolic operations is therefore not just an optimization but a critical strategy to move problems, such as the one treated in~\cite{santos2022decomposition}, from the realm of the intractable to the feasible.

 The constructive framework for systematically carrying out such symmetry reductions is developed comprehensively in many classical references, which provide a structured approach from representation-theoretic analysis to the explicit construction of symmetry-adapted bases.
 The cost of the classical implementation of this method grows quickly with the order of the group, a challenge that recent research has addressed with new algorithms that achieve polynomial-time guarantees for specific groups \cite{olver2025representations}.

 Although powerful computational algebra systems such as \texttt{GAP} \cite{gap} and \texttt{Magma} \cite{MR1484478}, along with the open-source \texttt{SageMath} environment \cite{sagemath}, provide extensive capabilities for computational group theory, none of them supplies a dedicated, end-to-end implementation of this specific workflow. Dedicated symmetry-reduction packages do exist, but each is confined to a single computational ecosystem and to one side of the symbolic--numerical divide: \texttt{GTPack} \cite{geilhufe2018gtpack} offers symbolic workflows for finite groups in \texttt{Mathematica}, \texttt{RepLAB} \cite{RepLAB} implements a numerical commutant-based approach for compact groups in \texttt{MATLAB}/\texttt{Octave}, and the representation-theory packages published most recently in this journal \cite{nazarov2012affine, feger2020lieart, fonseca2021groupmath} compute structural data for Lie algebras rather than acting on a user-supplied equivariant operator. Section~\ref{sec:related} places these tools and the \texttt{Maple} \texttt{GroupTheory} package side by side with \texttt{PySymmetry} on a common set of axes (Table~\ref{tab:tool_comparison}), and states where \texttt{PySymmetry} is the weaker choice.

Measured against that state of the art, this paper makes three contributions. First, \texttt{PySymmetry} carries a single constructive, character-theoretic method through two paths inside one open environment: an exact path, in which the projectors are built over an exact ring so that the resulting basis is free of round-off, and a numerical path interoperable with \texttt{NumPy}/\texttt{SciPy} and with sparse matrix formats. Rather than a domain-specific toolbox, it is a general-purpose workflow framework, and no package we are aware of offers both paths without leaving its host system. Second, this hybrid design delivers a tractability result rather than a speedup: the symbolic path produced the first complete analytical classification of the inverse problem for the nested cube central configuration in celestial mechanics, a computation that intermediate expression swell had previously placed out of reach. Third, on a $D_4$-equivariant non-Hermitian Schrödinger eigenvalue problem the numerical path attains a speedup of $16.93 \pm 1.45$ over a general-purpose dense eigensolver on the largest grid tested, with the gain growing as the grid is refined (Table~\ref{tab:performance_schrodinger}).

\texttt{PySymmetry} automates the symmetry reduction of
$G$-equivariant linear systems by providing a robust and reproducible implementation. The framework's core operation—generating a symmetry-adapted basis to transform an equivariant operator into its block-diagonal form—addresses both computational paradigms: in numerical contexts it replaces a single, expensive operation with a set of more manageable ones, and in symbolic contexts it acts as an enabler, preemptively mitigating expression swell by replacing a monolithic symbolic matrix with a collection of tractable blocks. Furthermore, this decomposition into independent subsystems inherently exposes a natural path for parallel computation, a powerful advantage that benefits both the numerical and symbolic approaches.

We demonstrate \texttt{PySymmetry}'s effectiveness and versatility through three comprehensive case studies: (1) a simple example from chemistry illustrating the GF method~\cite{WilsonFG, WilsonCourtneyCross} to calculate vibrational frequencies of molecules;
(2) a numerical benchmark on a $G$-equivariant, non-Hermitian Schrödinger operator, designed to demonstrate the significant performance advantages of symmetry reduction in computationally demanding scenarios; and (3) a symbolic investigation of the nested cube central configuration in celestial mechanics. By bridging abstract group theory with practical computational workflows, \texttt{PySymmetry} aims to serve as a robust and accessible tool for researchers across multiple disciplines.

The remainder of the paper follows the two paths that the framework unifies. Section~\ref{sec:related} places \texttt{PySymmetry} among the existing tools and states what none of them offers. Sections~\ref{sec:foundations} and~\ref{sec:reduction} set out the character-theoretic construction that the implementation follows, and Section~\ref{sec:framework} describes how that construction is realized both over an exact ring and over sparse numerical formats. Section~\ref{sec:applications} exercises the two paths on the three case studies, Section~\ref{sec:discussion} reads the results against the cost of the reduction, Section~\ref{sec:limitations} states where the framework does not reach, and Section~\ref{sec:conclusions} states what the two paths together deliver.

\section{Related Work — Software Landscape Across Computational Ecosystems}
\label{sec:related}

Symmetry-reduction capabilities appear in several computational environments, each shaped by the host system it was written for. We compare them here on a fixed set of axes: ecosystem, license and access, core approach, group scope, the route taken to block-diagonal form, ideal use cases, main advantage, and limitations. Table~\ref{tab:tool_comparison} summarizes the comparison. The paragraphs below supply the detail, and the last two state where \texttt{PySymmetry} is the weaker choice and why no cross-tool runtimes are reported.

\texttt{GTPack} \cite{geilhufe2018gtpack, hergert2018group} is a \texttt{Mathematica}-based suite providing advanced tools for solid-state physics and photonics, aimed in particular at classical problems in these fields, such as crystal field splitting and photonic band structures. Despite its specialized modules, \texttt{GTPack} rests on a robust, general-purpose mathematical foundation: it includes subpackages for abstract group theory and representation theory, with which the user can compute character tables, irreducible representations, and Clebsch--Gordan coefficients. Its modular architecture allows custom groups to be installed through multiplication tables, so that the software applies to any set of matrices forming a group, extending its utility beyond the predefined configurations.

\texttt{RepLAB} \cite{RepLAB} is designed for \texttt{MATLAB} or \texttt{Octave}, and so shares with \texttt{PySymmetry} the advantage of being usable on an open-source platform. It provides dedicated apparatus for optimization and convex programming problems, which distinguishes it from \texttt{PySymmetry}. While \texttt{RepLAB} follows \texttt{MATLAB} syntax, it shares several capabilities with our framework: building permutation groups from sets, using custom group constructors, working with tensor products and direct sums, and extracting irreducible subrepresentations. It further provides modules associating groups with their undirected graphs, allows computations with continuous groups, and can be integrated with the \texttt{GAP} system after specific configuration adjustments. Its route to block-diagonal form is a numerical decomposition of the commutant, which makes it stable and scalable for large-scale numerical work and, correspondingly, leaves it without an exact symbolic path.

The role of \texttt{Maple} in this landscape is a more delicate one. \texttt{Maple} ships a built-in \texttt{GroupTheory} package \cite{postma2014maple} covering the construction and structural analysis of finite groups: subgroups, conjugacy classes, homomorphisms, and character tables, backed by the system's exact symbolic arithmetic. That is precisely the input data our workflow consumes. What we were unable to find, either in the published literature or among the distributed packages, is a \texttt{Maple} counterpart to \texttt{GTPack} or \texttt{RepLAB}: a package that carries this data forward into the construction of a symmetry-adapted basis and the block-diagonalization of a user-supplied equivariant operator. In \texttt{Maple} that step is at present assembled by hand for each problem, which is why the system appears in Table~\ref{tab:tool_comparison} as a source of structural data rather than as a reduction pipeline.

A similar boundary separates \texttt{PySymmetry} from the representation-theory packages published most recently in this journal. \texttt{Affine.m} \cite{nazarov2012affine} computes weight multiplicities, branching rules, and tensor-product decompositions for finite-dimensional and affine Lie algebras; \texttt{LieART 2.0} \cite{feger2020lieart} covers tensor products and subalgebra branching for all classical and exceptional Lie algebras; \texttt{GroupMath} \cite{fonseca2021groupmath} performs calculations on semisimple Lie algebras and on permutation groups for particle-physics model building. All three are hosted in \texttt{Mathematica}, and all three compute structural data about groups and their representations. None of them takes a concrete $G$-equivariant operator and returns the change of basis that block-diagonalizes it, which is the operation \texttt{PySymmetry} automates.

All of these frameworks offer comprehensive tools for studying groups, yet they operate under distinct design philosophies. \texttt{GTPack} is aimed primarily at physicists and chemists dealing with applied problems, whereas \texttt{PySymmetry} and \texttt{RepLAB} address a broader audience and allow researchers to approach problems of a purely mathematical nature as well. \texttt{PySymmetry} occupies the position left open by the others: a constructive, character-theoretic method able to yield exact projector-based bases when symbolic derivations are required, while maintaining full interoperability with \texttt{NumPy}/\texttt{SciPy} for large-scale numerical problems. A major strength is its seamless interoperability with the Python libraries and its ability to invoke finite groups already present in \texttt{SageMath}, such as the dihedral and cyclic families. Its native integration with \texttt{GAP} gives access to specialized packages in turn.

\begin{table}[H]
\centering
\footnotesize
\renewcommand{\baselinestretch}{1}\selectfont

\caption{Comparative analysis of symmetry-reduction software packages, on a common set of axes. Entries are compiled from the cited publications and the package documentation \cite{geilhufe2018gtpack, hergert2018group, RepLAB, postma2014maple}. No cross-tool runtimes were measured, for the reason given in the text.}
\label{tab:tool_comparison}

\newcolumntype{L}{>{\raggedright\arraybackslash}X}

\setlength{\emergencystretch}{2em}

\begin{tabularx}{\textwidth}{@{}>{\raggedright\arraybackslash}p{0.135\textwidth} L L L L@{}}
\toprule
\textbf{Feature} & \textbf{PySymmetry} & \textbf{GTPack} & \textbf{RepLAB} & \textbf{Maple}\newline \textbf{GroupTheory} \\
\midrule
Ecosystem & Python / SageMath & Mathematica & MATLAB / GNU Octave & Maple \\
\addlinespace
License / access & Open-source (MIT) & Free package; requires proprietary Mathematica & Open-source (MPL 2.0), compatible with Octave & Built into the system; requires proprietary Maple \\
\addlinespace
Core approach & Hybrid (symbolic and numeric) & Predominantly symbolic & Strictly numeric & Symbolic; structural only \\
\addlinespace
Group scope & Finite groups & Finite groups & Compact groups & Finite groups \\
\addlinespace
Route to block-diagonal form & Character projection operators & Character projection operators & Numerical decomposition of the commutant & Not provided \\
\addlinespace
Ideal use cases & Hybrid equivariant systems, expression swell mitigation, symbolic and algebraic treatment of equivariance & Solid-state physics, photonics, analysis of symmetry-adapted functions & Invariant semidefinite programs (SDP), compact groups & Abstract structural group theory \\
\addlinespace
Main advantage & Full interoperability with the SageMath ecosystem, which includes NumPy, SciPy, GAP and PARI & Deep integration with Mathematica's ecosystem, databases and visual tools & Stability for large-scale numeric computation and integration with convex optimization (YALMIP) & Exact structural data available without any add-on \\
\addlinespace
Limitations & Restricted to finite groups; no domain-specific modules; requires a SageMath environment & Restricted to finite groups; requires a proprietary Mathematica license & Strictly numerical (lacks exact symbolic computation) & No construction of symmetry-adapted bases or block-diagonalization \\
\bottomrule
\end{tabularx}
\end{table}

Despite its strengths, \texttt{PySymmetry} has definite limitations when contrasted with the tools above, and the axes on which it falls short are worth stating plainly. Unlike its counterparts, it lacks built-in, domain-specific modules for applied physics or chemistry, functioning strictly as a general-purpose tool. For example, both \texttt{GTPack} and \texttt{RepLAB} include native routines for Clebsch--Gordan coefficients, and that functionality is currently absent here, although it and other specific features can be added in future versions. \texttt{RepLAB} treats compact groups, a class that properly contains the finite groups \texttt{PySymmetry} supports. For continuous symmetry, and for invariant semidefinite programs at scale, \texttt{RepLAB} is the appropriate choice and \texttt{PySymmetry} is not a substitute for it. \texttt{GTPack} ships curated crystallographic point- and space-group data that \texttt{PySymmetry} does not provide, so a solid-state workflow already resting on those tables is served more directly there. \texttt{Maple} and the \texttt{Mathematica} Lie-algebra packages remain stronger for the structural questions they were written to answer. Within its own scope, the cost of the \texttt{PySymmetry} workflow is dominated by the construction of the symmetry-adapted basis rather than by the reduced solve (Table~\ref{tab:performance_schrodinger}), which makes the method unattractive for a single small system. The investment pays off when the basis is reused across a parametric family, when the problem is large, or when the alternative is a symbolic computation that does not terminate.

Symmetry is also an active methodological line in computational structural mechanics, on a track parallel to the software one traced above. Zingoni surveys how group-theoretic decomposition is applied across solid and structural mechanics~\cite{Zingoni2009}, and Kaveh's monograph develops the analysis of structures from symmetry and regularity, including the graph-product constructions that generate the regular models such methods target~\cite{kaveh2013optimal}. A second strand attacks the step that precedes any reduction, namely finding the symmetry in the first place. Chen et al.\ detect cyclic symmetries in engineering structures automatically~\cite{chen2017cyclic}, and Zhang et al.\ recognize structural symmetry in planar structures with convolutional neural networks~\cite{zhang2022cnn}. That strand is complementary to the present work rather than competing with it, since \texttt{PySymmetry} takes the group as given, and it is the natural source for the automatic detection discussed in Section~\ref{sec:limitations}.

The comparison above is one of capability and design rather than of measured runtimes. A single benchmark shared across four host systems would confound the reduction method with three unrelated runtimes and linear-algebra backends, and a timing produced by the authors of one of the tools carries limited evidential weight; every entry in Table~\ref{tab:tool_comparison} is instead verifiable against the cited publications and documentation.

\section{Mathematical Foundations of Symmetry Reduction}
\label{sec:foundations}

The systematic exploitation of symmetry in linear problems hinges on the mathematical framework of group representation theory. This section outlines the core concepts, primarily following the constructive approach of Fässler and Stiefel \cite{stiefel2012group} and Serre \cite{serre1977linear}, which form the theoretical basis for the \texttt{PySymmetry} framework.

\subsection{Representations and Reducibility}
Let $G$ be a finite group and $V$ be a finite-dimensional vector space over $\mathbb{C}$. A \emph{linear representation} of $G$ is a homomorphism $\phi: G \rightarrow GL(V)$, where $GL(V)$ is the group of invertible linear transformations on $V$. A subspace $U \subseteq V$ is called \emph{$G$-invariant} if $\phi(g)u \in U$ for all $g \in G$ and $u \in U$. A representation $\phi$ is \emph{irreducible} if its only $G$-invariant subspaces are $\{0\}$ and $V$ itself; otherwise, it is \emph{reducible}.

A cornerstone of the theory for finite groups is \textbf{Maschke's Theorem}, which states that any reducible representation can be decomposed into a direct sum of irreducible subrepresentations. This property, known as complete reducibility, guarantees that any representation space $V$ can be written as:
\begin{equation}
    V = \bigoplus_{i=1}^{k} U_i,
\end{equation}
where each $U_i$ is an irreducible $G$-invariant subspace.

\subsection{Isotypic Decomposition and \texorpdfstring{$G$}{G}-Equivariant Operators}
The decomposition into irreducible subspaces is generally not unique. However, by grouping equivalent irreducible subrepresentations, we can obtain a unique canonical decomposition. The \emph{isotypic component} $V_j$ is the direct sum of all irreducible subspaces in $V$ that are equivalent (isomorphic as $G$-modules) to a specific irreducible representation (irrep) $\vartheta_j$. This yields the unique \emph{isotypic decomposition} of the space:
\begin{equation}
    V = \bigoplus_{j=1}^{N} V_j,
\end{equation}
where $N$ is the number of non-equivalent irreducible representations of $G$. If the irreducible representation $\vartheta_j$ has degree (dimension) $n_j$ and appears in $\phi$ with multiplicity $c_j$, then the isotypic component $V_j$ has dimension $n_j c_j$.

Our primary interest lies in linear problems of the form $Mx = b$ or $Mx = \lambda x$, where the operator $M: V \rightarrow V$ possesses the symmetry of the group. We say $M$ is \emph{$G$-equivariant} if it commutes with the group action given by the representation, i.e.:
\begin{equation}
    M\phi(g) = \phi(g)M \quad \forall g \in G.
\end{equation}
The central result that enables computational simplification is the \textbf{Fundamental Theorem of Group Representation Theory}, which states that a $G$-equivariant operator $M$ is block-diagonalized by a \emph{symmetry-adapted basis}. Specifically, the operator $M$ maps each isotypic component $V_j$ into itself ($M(V_j) \subseteq V_j$). Furthermore, when expressed in a symmetry-adapted basis, the matrix of $M$ decomposes into blocks, with $n_j$ identical blocks of size $c_j \times c_j$ for each isotypic component $V_j$.

\subsection{Algorithmic Construction of a Symmetry-Adapted Basis}
The core of \texttt{PySymmetry} is the algorithmic construction of this symmetry-adapted basis, which allows for the practical application of the Fundamental Theorem. The procedure follows these steps:

\subsubsection{Step A: Analysis using Characters and Conjugacy Classes}
The first step is to determine the structure of the representation $\phi$, i.e., to find the multiplicities $c_j$ of each irreducible representation $\vartheta_j$. This is achieved using the \emph{characters} of the representations. The character $\chi_\phi: G \rightarrow \mathbb{C}$ is defined as $\chi_\phi(g) = \mathrm{tr}(\phi(g))$. For finite groups, the characters of the irreducible representations form an orthonormal basis for the space of class functions. The multiplicity $c_j$ is computed via the inner product of their characters:
\begin{equation}\label{eq:char}
    c_j = \langle \chi_\phi, \chi_j \rangle = \frac{1}{|G|} \sum_{g \in G} \chi_\phi(g) \overline{\chi_j(g)}.
\end{equation}
While theoretically fundamental, a direct implementation of Eq.~\eqref{eq:char} is computationally expensive for groups of high order, as it requires summing over all $|G|$ elements. However, since characters are class functions (i.e., constant on conjugacy classes), this sum can be reformulated into a much more efficient computation over the $N_{\mathrm{classes}}$ distinct conjugacy classes of $G$:
\begin{equation} \label{eq:char-classes}
    c_j = \frac{1}{|G|} \sum_{k=1}^{N_{\mathrm{classes}}} |C_k| \cdot \chi_\phi(g_k) \overline{\chi_j(g_k)},
\end{equation}
where $C_k$ is the $k$-th conjugacy class, $|C_k|$ is its size, and $g_k \in C_k$ is any representative element. This reformulation, which is standard in computational group theory software \cite{Dixon1967}, reduces the complexity of this stage from $O(|G|)$ to $O(N_{\mathrm{classes}})$, a critical optimization for scalability.

\subsubsection{Step B: Projection Operators}
The key tools for constructing the basis are the \emph{projection operators}, built from the matrix entries of the irreducible representations. Given a fixed set of matrix representations $\{D_j(g)\}$ for the irreducible representation $\vartheta_j$, where $D_j(g) = [d^j_{k\ell}(g)]$, the transference operators are defined as:
\begin{equation}\label{eq:projection}
    P_{k\ell}^j = \frac{n_j}{|G|} \sum_{g \in G} d^j_{\ell k}(g^{-1}) \phi(g),
\end{equation}
where $n_j$ is the degree of $\vartheta_j$. As with the character inner product, a direct summation over the group is computationally expensive, a fact noted early in the history of computational group theory \cite{dixon1970computing}. Modern implementations, such as that in \texttt{PySymmetry}, avoid this by pre-calculating the representation matrices only for a minimal set of elements (the class representatives and their inverses) and applying the operator efficiently. The operator $P_{kk}^j$ is a projector onto a $c_j$-dimensional subspace of the isotypic component $V_j$, and the operator $P_{k\ell}^j$ ($k \neq \ell$) is an isomorphism between two irreducible subrepresentations of $V_j$.

\subsubsection{Step C: Basis Generation and Ordering}
The algorithm implemented in \texttt{PySymmetry} generates the basis for each isotypic component $V_j$ as follows:
\begin{enumerate}
    \item A basis for the image of the projector $P_{11}^j$ is computed. Since the rank of this operator is the multiplicity $c_j$, this step yields $c_j$ linearly independent vectors: $\{v_1^1, v_1^2, \dots, v_1^{c_j}\}$. Each vector belongs to a different irreducible subspace within $V_j$.
    \item The remaining basis vectors are generated by applying the transference operators:
    \begin{equation}
        v_k^i = P_{1k}^j v_1^i \quad \text{for } k=2, \dots, n_j \text{ and } i=1, \dots, c_j.
    \end{equation}
    The set $\{v_1^i, v_2^i, \dots, v_{n_j}^i\}$ forms a basis for the $i$-th irreducible subspace $V_j^i$.
    \item The final crucial step is the ordering of the basis vectors. As illustrated in \cite{stiefel2012group, santos2022decomposition}, ordering the basis by rows (i.e., grouping vectors from the same irreducible subspace $V_j^i$ together) results in a matrix for $\phi(g)$ that is block-diagonal. However, to achieve the desired simplification of the equivariant operator $M$, the basis must be ordered by columns. This \emph{symmetry-adapted basis} is constructed by grouping vectors that correspond to the same column index $k$ across all the different irreducible subspaces: $\{\{v_1^1, \dots, v_1^{c_j}\}, \{v_2^1, \dots, v_2^{c_j}\}, \dots\}$. This reordering yields the block-diagonal form for $M$ predicted by the Fundamental Theorem.
\end{enumerate}
This systematic procedure provides a direct computational path from the definition of a group and its representation to the simplification of a corresponding equivariant linear system.

\section{Reducing Representations into Their Irreducible Components}
\label{sec:reduction}
\texttt{PySymmetry} is also a general-purpose tool for working with representations. Some of the classical algorithms presented in \cite{stiefel2012group, serre1977linear} are implemented in \texttt{PySymmetry}. Among its capabilities, we can list the following:
\begin{enumerate}
    \item Defines a representation by giving the representative matrices of the group generators.
    \item Finds the irreducible representations of a group and gives a complete set of representative matrices for each of them.
    \item Calculates the isotypic basis associated with a representation.
    \item Calculates a basis that decomposes a representation into its irreducible components.
    \item Calculates the symmetry-adapted basis for an operator that is equivariant under the representation.
    \item Calculates the multiplicities of the irreducible subrepresentations of a representation.
    \item Based on representation degrees and multiplicities, gives the size of the block decomposition for an equivariant operator.
    \item Includes many relevant tests, for example,
     whether a map is indeed a representation of the defined group
     and whether an operator is equivariant under a representation.
\end{enumerate}

Some features and methods of \texttt{PySymmetry} are described and discussed below. The reader should note that this is not an exhaustive list.

\texttt{PySymmetry} allows the user to define representations over many of the fields available in the \texttt{SageMath} environment. This makes it possible to work with numeric fields as well as symbolic fields. The main properties of a representation are available through many methods. 

One can create a group $G$ by providing a set of generators as permutation elements. Given the group, one can generate the complete set of irreducible representations of $G$ with the method \verb|irreducible_representations|. One can also create a representation in matrix form from this group to any general linear group $GL_n(\mathbb{C})$ using the method \verb|representation|.

It is also possible to verify whether the defined map satisfies the requirements of a representation using the method \verb|is_representation|. Once the representation is defined, it is possible to verify whether a given operator is equivariant with respect to the defined representation using the method \verb|is_equivariant_to|.
Based on the properties of the representation, one can also preview the sizes of the blocks in the block decomposition of the equivariant operator via the method \verb|quick_block_prevision|.
If one only needs the restriction of the equivariant operator to a specific subspace, one can use the method \verb|base_change_equivariant_reduced| to get the list of invariant subspaces together with the method \verb|get_block| to get the actual expression of the restriction of the equivariant operator to the subspace.

 \texttt{PySymmetry} can reduce a representation into its irreducible components, and the method \verb|base_to_irreducibles| returns the explicit matrix that makes the change of basis.

Because of its relevance, the regular representation is built into \texttt{PySymmetry} and can be called through the method \verb|regular_representation|.

Permutation representations are also available in \texttt{PySymmetry} and can be invoked via the method \verb|natural_representation|.

\texttt{PySymmetry} is also equipped with methods that allow one to define new representations from others already defined. For example, one can form the direct sum and the tensor product of representations through the methods \verb|direct_sum| and \verb|tensor_product|, respectively.
Other specific features of \texttt{PySymmetry} are described in the package documentation.

\section{The 	\texttt{PySymmetry} Framework: Architecture and Implementation}
\label{sec:framework}

\texttt{PySymmetry} is designed to be a practical and accessible tool for applying the powerful methods of group representation theory. This section details its architecture, the implementation of its core algorithms, and the workflow for symmetry reduction.

\subsection{Design Philosophy: Bridging Symbolic Algebra and Numerical Computation}
The \texttt{PySymmetry} framework is strategically built upon the \texttt{SageMath} environment. This choice was motivated by the unique ability of \texttt{SageMath} to bridge the gap between high-level symbolic algebra and high-performance numerical computation. For many problems in group theory, the
ability to perform exact calculations on symbolic rings is crucial for correctness. \texttt{PySymmetry} leverages \texttt{Sage}'s powerful symbolic engine and its direct interface to specialized computer algebra systems like \texttt{GAP} \cite{gap} for these tasks, while incorporating advanced algorithmic optimizations to ensure computational feasibility even for large-order groups.

Simultaneously, the framework retains full access to the extensive Python scientific ecosystem, including libraries like \texttt{NumPy} and \texttt{SciPy}. This hybrid approach is a key advantage, allowing \texttt{PySymmetry} to be applied to problems requiring exact symbolic results and to large-scale numerical problems.

\subsection{Core Architecture and Data Structures}
The framework is designed with an object-oriented architecture to encapsulate the core mathematical concepts of representation theory. The main classes are:
\begin{itemize}
    \item \textbf{\texttt{FiniteGroup}}: This class serves as the foundational object, inheriting from \texttt{SageMath}'s \texttt{PermutationGroup\_generic}. It provides a consistent interface for defining a finite group and accessing its intrinsic properties, such as its elements, its order, and its irreducible representations via an interface to \texttt{GAP}.
    \item \textbf{\texttt{MapRepresentation}}: This class implements the concept of a linear representation $\phi$. It acts as a map from the \texttt{FiniteGroup} object to a matrix group, effectively defining the vector space where the representation takes place.
    \item \textbf{\texttt{IsotypicBase}}: This class is the primary data structure for the output. It holds the final change-of-basis matrix $P$ and stores the metadata of the decomposition, including the degree ($n_j$) and multiplicity ($c_j$) of each block, providing a clear overview of the decomposition structure.
\end{itemize}

A key feature of the architecture is its flexibility in matrix representation. As demonstrated in the associated \texttt{pde.py} module, numerical representations using formats like \verb|scipy.sparse.csr_matrix| can be employed, which is critical for achieving computational efficiency in large, sparse linear systems. To distinguish these numerically optimized routines from their symbolic counterparts, the core user-facing classes and functions intended for performance benchmarks, such as \texttt{nFiniteGroup} and \verb|nbase_change_reduction|, are prefixed with \texttt{n}. This design choice ensures a clear separation between the symbolic and numerical workflows within the framework.

\subsection{The Algorithmic Workflow for Symmetry Reduction}
\label{sec:workflow}
The main user-facing function, \texttt{base\_equivariant\_to\_blocks}, orchestrates the entire symmetry reduction pipeline, directly implementing the algorithmic steps derived from the theoretical framework.
\begin{enumerate}
    \item \textbf{Analysis Stage:} The workflow begins by determining the structure of the input representation $\phi$. A naive implementation would calculate the multiplicities $c_j$ by summing over all $|G|$ group elements. However, to enable the analysis of systems with large symmetry groups, \texttt{PySymmetry} implements a crucial algorithmic optimization. Taking advantage of the fact that characters are class functions, the \texttt{inner\_product} method computes the multiplicities via a much more efficient sum over the group's conjugacy classes. This reduces the computational complexity of the analysis stage from being proportional to the group order, $|G|$, to being proportional to the number of conjugacy classes, $N_{\mathrm{classes}}$, which is often orders of magnitude smaller.

    \item \textbf{Basis Generation via Projection:} The core of the framework is the construction of the symmetry-adapted basis. This process, encapsulated within the \texttt{isotypic\_base} method, implements an efficient procedure. The method first computes the matrix representation of the projection operator $P_{11}^j$. It then calculates a basis for the image (or column space) of this matrix using standard linear algebra routines. For symbolic computations, this is achieved exactly by finding the pivot columns via Gaussian elimination within \texttt{SageMath}. Since the rank of the projector is known to be the multiplicity $c_j$, this step directly yields a set of $c_j$ linearly independent starting vectors, $\{v_1^1, v_1^2, \dots, v_1^{c_j}\}$. The full basis for the isotypic component is then generated by applying the transference operators to this well-defined initial set.

    \item \textbf{Output and Block-Diagonalization:} The generated basis vectors are systematically reordered by their column index ($k$) to form the final symmetry-adapted basis matrix $P$, returned within an \texttt{IsotypicBase} object. The simplified, block-diagonal form of the original equivariant operator $M$ is then obtained via the similarity transformation $\tilde{M} = P^{-1}MP$.
\end{enumerate}

\subsection{Illustrative Workflow}
To make the algorithmic process concrete, we provide a detailed, step-by-step example in \ref{app:d4_example}. The example demonstrates the complete workflow of the framework on the symmetry reduction of a $D_4$-equivariant operator. This walkthrough showcases the high-level API of \texttt{PySymmetry}, which abstracts away the complex underlying calculations, allowing the user to focus on the problem's definition and the analysis of the simplified result.
\section{Applications and Case Studies}
\label{sec:applications}

This section demonstrates the practical value and versatility of the \texttt{PySymmetry} framework through three distinct case studies. 
First, we give a simple example of how \texttt{PySymmetry} can be used to carry out the calculations of an important method from chemistry.
Second, we establish the framework's correctness and computational efficiency on a canonical numerical problem. Third, we showcase its power as an enabling technology for solving a previously challenging problem in symbolic celestial mechanics.

\subsection{Computational Environment}
\label{sec:environment}

All timings reported in this paper were measured on a single machine, whose specification is given here so that they can be placed in context and reproduced. The accuracy measurements reported at the end of Section~\ref{sec:schrodinger} were produced on a different machine and carry no timing.

\begin{itemize}
    \item \textbf{Processor:} 12th generation Intel Core i7-1255U, 10 cores and 12 threads, base clock 1.70 GHz.
    \item \textbf{Memory:} 16.0 GB installed, 15.7 GB usable.
    \item \textbf{Operating system:} Ubuntu 22.04.4 LTS under the Windows Subsystem for Linux 2, hosted on Windows 11 Home Single Language, version 24H2.
    \item \textbf{Software stack:} \texttt{SageMath} 10.4, which integrates Python 3.11, together with \texttt{NumPy} 1.26 and \texttt{SciPy} 1.14. The block eigenvalue calculations were distributed over the available CPU threads with Python's \texttt{multiprocessing} library.
\end{itemize}

\subsection{An Illustrative Application in Chemistry: The GF Method}
\label{sec:gf}
The GF method of Wilson \cite{WilsonFG, WilsonCourtneyCross} is a classical method for studying molecular vibrations.
To calculate the vibrational frequencies, one is interested in finding the eigenvalues of the product of matrices $\textbf{F}$ and $\textbf{G}$, where $\textbf{F}$ is the force-constant matrix, and $\textbf{G}$ represents the inverse kinetic energy matrix.
The $\mathbf{F}$ matrix is derived from the molecular potential energy, while the $\mathbf{G}$ matrix depends on the atomic masses and the molecular geometry.
When the molecule has symmetry, the secular equation
\begin{equation}
\det(\textbf{FG}-\lambda I) = 0
\label{secular_equation}
\end{equation}
factors into smaller equations, because $\textbf{FG}$ can be block-diagonalized by transforming to a symmetry-adapted basis. This simplification is possible because both the $\mathbf{F}$ and $\mathbf{G}$ operators commute with the symmetry operations of the molecule's point group.

As an illustrative example, consider the water molecule $\mathrm{H_2O}$. This is a planar molecule. We place the oxygen atom at the origin of a Euclidean space with axes $x$, $y$, and $z$, and the hydrogen atoms in the $zx$-plane. In this case, $\textbf{F}$ and $\textbf{G}$ are symmetric and have the form
\begin{align}
 \mathbf{F} = \begin{pmatrix}
f_{11} & f_{12} & 0 \\
f_{12} & f_{11} & 0 \\
0 & 0 & f_{33}
\end{pmatrix},\,\mathbf{G} = \begin{pmatrix}
g_{11} & g_{12} & g_{13} \\
g_{12} & g_{11} & g_{13} \\
g_{13} & g_{13} & g_{33}
\end{pmatrix}.   
\end{align}

The water molecule has the symmetry of the group $C_{2v}$.
In the standard notation of chemistry, the four elements of $C_{2v}$ are $E$ (the identity), $C_2$
(the rotation by $\pi$ about the $z$-axis), $\sigma_v(xz)$ (the reflection through the $zx$-plane) and
$\sigma'_v(yz)$ (the reflection through the $zy$-plane).

We can think of $C_{2v}$ as a permutation group simply by numbering the atoms, for example (H, H, O) = (1, 2, 3), and
associating each group element with the permutation it causes on the atoms, as follows:
$E \to (1)(2)(3)$,
$C_2 \to (12)(3)$,
$\sigma_v(xz) \to (1)(2)(3)$ and
$\sigma'_v(yz) \to (12)(3)$.
We can then define a permutation representation that encodes the change in the atoms' positions under the elements of $C_{2v}$ as follows:
\begin{align}
(1)(2)(3) \longrightarrow \begin{pmatrix} 1 & 0 & 0 \\ 0 & 1 & 0 \\ 0 & 0 & 1 \end{pmatrix}, \,
(12)(3) \longrightarrow \begin{pmatrix} 0 & 1 & 0 \\ 1 & 0 & 0 \\ 0 & 0 & 1 \end{pmatrix} 
\label{representation_of_C2V}.
\end{align}
Using \texttt{PySymmetry}, we can decompose the secular equation \eqref{secular_equation} by finding the symmetry-adapted basis with the method \verb|base_equivariant_to_blocks|, obtaining
\begin{equation}
 \beta = \left(\begin{array}{rrr}
2 & 0 & 2 \\
2 & 0 & -2 \\
0 & 4 & 0
\end{array}\right).   
\end{equation}
In the new basis $\beta$, the matrix $\textbf{FG}$ has the block-diagonal form shown below, which simplifies the calculation of the vibrational frequencies.

\begin{align}
\beta^{-1}\mathbf{FG}\beta =
\left(
\begin{array}{c|c}
\begin{array}{rr}
{\left(f_{11} + f_{12}\right)} {\left(g_{11} + g_{12}\right)} & 2 \, {\left(f_{11} + f_{12}\right)} g_{13} \\
f_{33} g_{13} & f_{33} g_{33}
\end{array}
& \begin{array}{r}
0 \\[6pt] 0
\end{array} \\ \hline
\begin{array}{rr}
0 & 0
\end{array}
& {\left(f_{11} - f_{12}\right)} {\left(g_{11} - g_{12}\right)}
\end{array}
\right).
\end{align}

We refer the reader to \ref{app:GF_example} for more details of this example in the \texttt{PySymmetry} framework.

\subsection{Validation and Performance on a Canonical PDE Problem}
\label{sec:schrodinger}

To demonstrate the framework's capabilities, we first validate its correctness on a well-understood Hermitian system and then benchmark its performance on a more demanding, non-Hermitian yet fully $G$-equivariant case.

\begin{figure}[H]
    \centering
    \includegraphics[width=0.8\textwidth]{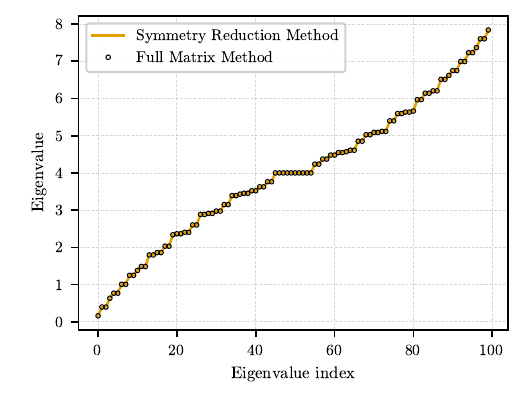}
    \caption{Correctness validation: sorted eigenvalue spectrum of the five-point Laplacian stencil on a $10 \times 10$ grid with unit spacing ($N=100$). The orange line is the symmetry-reduced spectrum returned by the \texttt{PySymmetry} workflow and the black circles are the full-matrix reference. The two coincide at every index, including the tenfold degeneracy at eigenvalue $4$. The agreement is quantified in the text.}
    \label{fig:eigenvalue_comparison}
\end{figure}

First, for correctness validation, we use the eigenvalue problem for the 2D discretized Laplacian operator, obtained on a uniform square grid using the classical two-dimensional finite-difference scheme (5-point stencil). The resulting matrix is Hermitian and  preserves the $D_4$
 symmetry. We computed the spectrum for a $10 \times 10$ system ($N=100$) and compared two numerical approaches. The reference spectrum was obtained by diagonalizing the full matrix using \texttt{numpy.linalg.eigvalsh}, a solver highly optimized for Hermitian matrices. This was compared with the spectrum from the \texttt{PySymmetry} workflow, for which the eigenvalues of the decomposed blocks were calculated using the general-purpose \texttt{numpy.linalg.eigvals} function to ensure numerical robustness. The sorted spectra agree to machine precision (Figure~\ref{fig:eigenvalue_comparison}), with a maximum absolute deviation of $1.3 \times 10^{-14}$ and a maximum relative deviation of $3.2 \times 10^{-14}$ over the $100$ eigenvalues, which validates the correctness of the \texttt{PySymmetry} framework's implementation.

Having established the framework's fidelity, we turn to a more rigorous benchmark to showcase its computational advantage: the eigenvalue problem for the 2D Schrödinger equation with a complex potential \cite{moiseyev2011non}, $H = -\nabla^2 + iV(x,y)$. This transition is critical for demonstrating the broader utility of \texttt{PySymmetry}. While performance gains against specialized Hermitian eigensolvers like \texttt{numpy.linalg.eigvalsh} are often modest, the true computational advantage of symmetry reduction is realized in non-Hermitian problems. Such systems require general-purpose solvers, like \texttt{numpy.linalg.eigvals}, which are inherently more computationally expensive. To ensure the operator remains perfectly \textbf{$D_4$-equivariant}, we chose a potential that is itself $D_4$-invariant, $V(x,y) = V_0(x^2 + y^2)$. This provides an ideal test case where the theory applies exactly, allowing for a direct and compelling performance comparison. For our experiments, we set the potential amplitude to $V_0=10.0$. Each benchmark was executed five times on the machine of Section~\ref{sec:environment}, and the mean and standard deviation of the timings were recorded.

The performance metrics for this problem are summarized in Table~\ref{tab:performance_schrodinger}. For the largest grid tested ($90\times 90$), the full-matrix method required an average of 704.55 seconds (nearly 12 minutes), while the \texttt{PySymmetry} approach completed the calculation in 41.88 seconds on average. This is a performance increase of nearly \textbf{17 times}.

\begin{table}[H]
\centering
\caption{Performance results (mean $\pm$ std over 5 runs) for the $D_4$-equivariant Schrödinger eigenvalue problem. Abbreviations: $T_{p}$ (preprocessing time), $T_{b}$ (block eigenvalue calculation time), $T_{s}$ (total time for the symmetry method), and $T_{f}$ (full-matrix method time). All times are in seconds.}
\label{tab:performance_schrodinger}

\newcolumntype{R}{>{\raggedleft\arraybackslash}X}

\small
\setlength{\tabcolsep}{3pt}
\begin{tabularx}{\textwidth}{c R R R R R} 

\toprule
\textbf{Grid} &  \textbf{$T_{p}$ (s)} & \textbf{$T_{b}$ (s)} & \textbf{$T_{s}$ (s)} & \textbf{$T_{f}$ (s)} & \textbf{Speedup} \\
\midrule

 $60\times 60$ &   5.10 $\pm$ 0.35 &  1.23 $\pm$ 0.17 &  6.33 $\pm$ 0.50 &  46.98 $\pm$ 2.83 &  7.47 $\pm$ 0.77 \\
 $70\times 70$ &   9.33 $\pm$ 0.56 &  2.76 $\pm$ 0.10 & 12.09 $\pm$ 0.52 & 124.37 $\pm$ 7.15 & 10.31 $\pm$ 0.79 \\
 $80\times 80$ &  15.79 $\pm$ 0.63 &  5.36 $\pm$ 0.14 & 21.16 $\pm$ 0.58 & 301.29 $\pm$ 3.51 & 14.25 $\pm$ 0.49 \\
$90\times 90$ &  30.98 $\pm$ 3.31 & 10.90 $\pm$ 0.81 & 41.88 $\pm$ 3.54 & 704.55 $\pm$ 28.05 & 16.93 $\pm$ 1.45\\
\bottomrule
\end{tabularx}
\end{table}

\begin{figure}[H]
    \centering
    \includegraphics[width=0.8\textwidth]{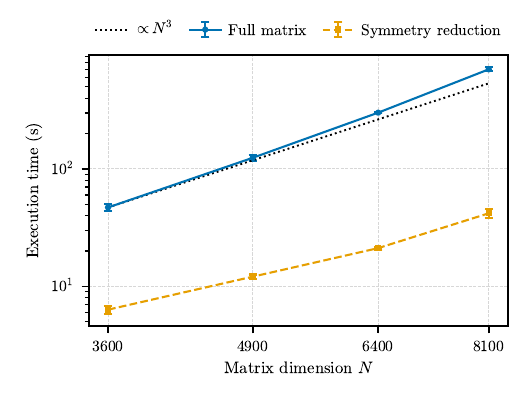}
    \caption{Execution time for the $D_4$-equivariant Schrödinger eigenvalue problem at four matrix dimensions, in log--log axes. The upper curve is the full-matrix method (\texttt{numpy.linalg.eigvals}); the lower one is the \texttt{PySymmetry} workflow timed end to end, the symmetry-adapted basis and the block solves together, written $T_{s}$ in Table~\ref{tab:performance_schrodinger}. Markers are means over five runs, the error bars extend one standard deviation either side of the mean, and the dotted line is a reference slope $\propto N^{3}$. The gap between the two curves widens with $N$.}
    \label{fig:performance_schrodinger}
\end{figure}

The two costs diverge steadily with $N$ (Figure~\ref{fig:performance_schrodinger}). The full-matrix runtime tracks the $\propto N^{3}$ reference, while the reduced path stays well below it. The reduction does not change the asymptotic order, since solving $k$ blocks of size $N/k$ still costs $O(N^{3})$, but it divides the constant by roughly $k^{2}$, and that is what the widening gap measures.

\begin{figure}[H]
    \centering
    \includegraphics[width=0.8\textwidth]{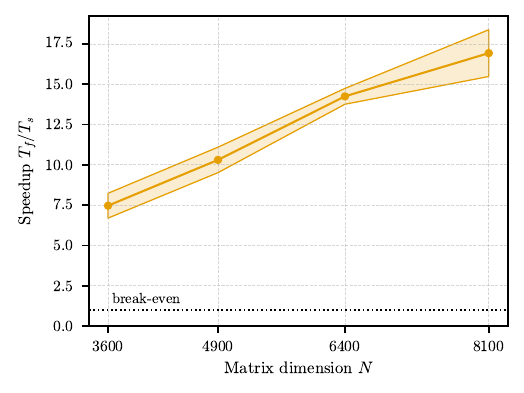}
    \caption{Speedup of the \texttt{PySymmetry} workflow over the full-matrix method for the $D_4$-equivariant Schrödinger problem, on the four grids and five runs of Figure~\ref{fig:performance_schrodinger}. The plotted ratio is the full-matrix time divided by the total time of the reduced path, $T_{f}/T_{s}$ in the notation of Table~\ref{tab:performance_schrodinger}, and the outlined band extends one standard deviation either side of the mean ratio. The dotted line marks break-even, where the reduction would cost as much as it saves. The speedup rises monotonically over the range tested.}
    \label{fig:speedup_schrodinger}
\end{figure}

The speedup grows with problem size (Figure~\ref{fig:speedup_schrodinger}), reaching an average of nearly \textbf{17$\times$} for the $8100 \times 8100$ matrix. This follows from the exact applicability of representation theory, which decomposes the equivariant operator into truly independent, smaller subproblems. The gain per unit of $N$ nevertheless falls by $40\%$ between the last two grid pairs, from $2.61$ to $1.59$ per $1000$ in $N$, so these four points do not exclude an eventual plateau. They do not establish one either, the shortfall against the trend set by the earlier grids being $1.2$ standard deviations.

The construction raises two questions of numerical stability, one for each path. On the symbolic path no round-off enters, because the basis is built over an exact ring. On the numerical path the quantity that matters is the conditioning of each isotypic basis $V_i$, since the blocks are extracted through a pseudoinverse. For a permutation representation that conditioning is fixed by the orbit structure rather than by the size of the grid. Every column of $V_i$ is a sum over one orbit of grid points, and distinct orbits have disjoint supports, so $V_i^{\mathsf{T}}V_i$ is diagonal. Then $\kappa_2(V_i)$ is the square root of the ratio between the largest and the smallest stabilizer order among the orbits that the irreducible representation supports. Measured from $N=100$ to $N=8100$, the values are exactly $1$ and $\sqrt{2}$ and match that expression at every grid, while a grid carrying a point fixed by the whole group raises the trivial component to $2\sqrt{2}$, again as predicted. The reconstructed eigenpairs are backward stable, with a relative residual of $2\times10^{-14}$ at $N=8100$, a hundred times the double-precision unit round-off, and the assembled spectrum agrees with the dense reference to a relative $10^{-11}$. The argument assumes a permutation representation with integer projector entries, and Section~\ref{sec:limitations} states the boundary. These accuracy figures were measured on a machine other than the one of Section~\ref{sec:environment}, which does not affect them because none of them is a timing.

\subsection{Main Case Study: Solving a Symbolic Problem in Celestial Mechanics}
\label{sec:celestial}
The inverse problem for central configurations represents a challenge in celestial mechanics, primarily due to the complexity of its governing symbolic equations. While simpler cases have been analyzed, others remain difficult to address because the size of the expressions and equations involved makes them hard to study analytically, or even computationally in a symbolic way.

As an example of a significant step, we address the
nested cube problem. Its symmetry group is $O_h$, a group of order 48.
The equations that define the central configuration in this case are given by a matrix $S(c,t)$ of size $48\times 16$. This matrix depends on the variables $t$ and $c$, and one wishes to discuss its null space as a function of these variables. This requires a change of basis in a symbolic matrix, followed by the study of the resulting equations.

Directly analyzing this matrix is computationally prohibitive due to intermediate expression swell, a well-known phenomenon where the symbolic complexity of intermediate calculations, such as computing determinants of large minors, grows factorially and exhausts the resources of modern computer algebra systems \cite{vonZurGathen2013modern}.

The system matrix $S(c,t)$ satisfies the general equivariance condition $S\theta(g) = (\theta \otimes \rho)(g)S$, where $\rho\colon G\to O(3)$ is the action of the symmetry group $G$ on $\mathbb{R}^3$ and $\theta$ is the permutation representation that $\rho$ induces on the bodies.

\texttt{PySymmetry} was the key to overcoming this barrier. Its framework is explicitly designed to handle this kind of problem. By applying the workflow described in Section~\ref{sec:workflow}, the intractable $48 \times 16$ symbolic matrix was block-diagonalized into a set of small, independent blocks (two $2 \times 2$ blocks and two sets of identical $4 \times 2$ blocks). This reduction was the critical enabling step that made it computationally feasible to analyze the system's null space, leading directly to the first complete symbolic classification of masses for the nested cube central configuration, a new scientific result reported in \cite{santos2022decomposition}.
Table~\ref{Order_of_groups_and_matrices} lists the representations decomposed in the same work, which also treats nested tetrahedra and nested octahedra.

\begin{table}[H]

\begin{center}

\begin{tabular}{l|c|c|c|c|c}
& \text{Group} & \text{Group order} & \text{Size of $S$} & \text{Degree of $\theta$} & \text{Degree of $\theta\otimes\rho$} \\[2pt]\hline
\text{Tetrahedron} & $S_4$ & $24$ &  $24\times 8$ & $8$ & $24$ \\\hline
 \text{Octahedron} & $O_h$ & $48$ &   $36\times 12$ & $12$ & $36$    \\\hline
\text{Cube} & $O_h$ & $48$ & $48\times 16$ & $16$  & $48$ \\ \hline
\end{tabular}
\end{center}
\caption{Representations decomposed with \texttt{PySymmetry} in \cite{santos2022decomposition} for central configurations of two nested regular polyhedra. The Cube row is the nested cube problem treated in the text. For $N$ bodies, $\rho\colon G\to O(3)$ is the action of the symmetry group $G$ on $\mathbb{R}^3$, and $\theta$ is the permutation representation that $\rho$ induces on the $N$ bodies, so $\deg\theta=N$ and $\deg(\theta\otimes\rho)=3N$. The matrix $S$ of the central configuration equations has size $3N\times N$ and satisfies $S\theta(g)=(\theta\otimes\rho)(g)S$ for all $g\in G$.}
\label{Order_of_groups_and_matrices}
\end{table}

\texttt{PySymmetry} decomposes representations of degree higher than those in Table~\ref{Order_of_groups_and_matrices}. What limits the approach is the symbolic analysis of the blocks it returns, which remains demanding even after the decomposition.

\section{Discussion}
\label{sec:discussion}

The results presented in this paper demonstrate that \texttt{PySymmetry} successfully provides a practical and automated workflow for the symmetry reduction of $G$-equivariant linear systems. The framework's effectiveness was validated across distinct computational paradigms. In high-performance numerical problems, the Schrödinger equation benchmark achieved a satisfactory average speedup of nearly $17\times$ against standard general-purpose solvers. Even more critically, in symbolic contexts, the tool enabled the complete analytical solution of a problem in celestial mechanics previously considered challenging due to intermediate expression swell.

A key insight from the performance benchmark is that the computational cost of the \texttt{PySymmetry} workflow is dominated by the initial preprocessing stage—the construction of the symmetry-adapted basis. This upfront investment is rewarded by the subsequent, rapid solution of the block-diagonalized system. Crucially, this basis matrix depends only on the system's symmetry group and representation, not on the specific numerical values of the operator. This implies that for common scientific workflows, such as parametric studies or time-dependent simulations where the equivariant operator changes but the underlying symmetry does not, the expensive preprocessing step needs to be performed only once. The cost can then be amortized over numerous subsequent calculations, leading to an even more dramatic effective speedup. Furthermore, the block-diagonal structure generated by \texttt{PySymmetry} inherently exposes an opportunity for efficient parallel processing on independent tasks, providing a clear pathway to substantial performance scaling on modern hardware.

Our approach complements existing toolkits in distinct computational ecosystems. While \texttt{GTPack} provides symbolic workflows in \texttt{Mathematica} and \texttt{RepLAB} offers a mature numerical path in \texttt{MATLAB}/\texttt{Octave} for compact groups, \texttt{PySymmetry} delivers a hybrid symbolic–numeric method in the open Python/\texttt{Sage} environment. Table~\ref{tab:tool_comparison} sets the four environments against a common set of axes. The value of exploiting symmetry in high-performance computing is also underscored by state-of-the-art research, such as the \texttt{Symtensor} library by Gao et al.~\cite{Gao2023}.

As a complement to these specialized packages, \texttt{PySymmetry} is a general-purpose, accessible tool that bridges abstract group theory with practical computation. The finite groups it supports form a proper subclass of the compact groups treated by \texttt{RepLAB}, and its contribution lies elsewhere, in one character-theoretic construction that serves an exact symbolic path and a sparse numerical path inside a single open environment. That combination is what fills a distinct and important gap in the scientific software ecosystem. Where a problem calls for continuous symmetry, or for tabulated crystallographic data, the packages discussed in Section~\ref{sec:related} remain the better choice.

\section{Limitations and Future Work}
\label{sec:limitations}

The cost profile that makes the reduction attractive for repeated solves makes it unattractive for a single one. At the largest grid tested, the construction of the symmetry-adapted basis accounts for $30.98$ of the $41.88$ seconds spent by the reduced path (Table~\ref{tab:performance_schrodinger}), and the measured gain falls from $16.93$ to $7.47$ as the grid contracts from $90\times 90$ to $60\times 60$. A user who solves one small system once will not recover that preprocessing, and we have not located the size below which the reduction stops paying for itself. The amortization argument of Section~\ref{sec:discussion} is what carries the method for larger workflows, and it assumes a setting in which the operator changes while the symmetry does not.

The framework is restricted to finite groups. Its constructions run through character tables and permutation groups, so continuous and compact symmetry lies outside them, and Section~\ref{sec:related} identifies where other packages serve those problems better. The group must also be supplied by the user, as a set of generators together with the representation acting on the discretized problem. \texttt{PySymmetry} does not detect the symmetry of an operator it is given, and automatic detection of the invariance group from the operator or from the geometry that produced it~\cite{chen2017cyclic,zhang2022cnn} is the extension we consider most useful, since it would remove the step at which a user unfamiliar with representation theory is most likely to stop.

The reduction is exact only for exactly equivariant operators, which is why the benchmark of Section~\ref{sec:schrodinger} uses a $D_4$-invariant potential. Where the symmetry is approximate, as with a perturbed potential, a lattice defect, or a boundary condition that breaks the group, the transformed operator retains couplings between blocks whose magnitude is set by the perturbation. The decomposition then becomes an approximation, and the framework provides no estimate of the resulting error.

Two implementation choices bound the problem sizes reachable today. The symbolic path returns the change of basis as a dense matrix over an exact ring, so its memory grows with the square of the degree of the representation, and the blocks it produces can still resist analysis, as the largest representations of Section~\ref{sec:celestial} show. The numerical path stores the basis in compressed sparse form, which is what makes the grids of Section~\ref{sec:schrodinger} feasible, but that basis is assembled from pivot columns of the projectors and is not orthonormalized, so its conditioning is inherited from the representation rather than controlled. That conditioning is bounded by the orbit structure and measured in Section~\ref{sec:schrodinger} for the permutation representations used there. Outside that case the bound does not apply, and the conditioning of the basis is what decides whether the reduced spectrum can be trusted.

The parallelism that the block structure exposes is realized today only within a single machine. Independent blocks are distributed over local processes, and the framework carries no distributed-memory or accelerator backend, so the scaling anticipated in Section~\ref{sec:discussion} is a structural property of the decomposition and has not been measured beyond one node. Supplying such a backend, and measuring what it delivers, is the natural continuation of the performance study reported here.

\section{Conclusions}
\label{sec:conclusions}

The framework reduces a $G$-equivariant linear system to its isotypic blocks through a single character-theoretic construction, realized once over an exact ring and once over sparse numerical formats. Of the three case studies, two measure what that construction delivers. A central-configuration problem in celestial mechanics, which intermediate expression swell had kept out of reach, was reduced to blocks small enough for a complete symbolic classification of masses, and a non-Hermitian Schr\"odinger operator was solved about seventeen times faster than by a general-purpose dense eigensolver at the largest grid tested, on a reduced basis whose conditioning is fixed by the orbits of the group rather than by the size of the grid.

Both paths live in one open environment, so a user can move from the exact treatment of a small system to the numerical treatment of a large one without reformulating the problem. The framework should be most useful where the operator changes while the symmetry does not, and the two extensions that would widen its reach, a distributed backend and automatic detection of the invariance group, are the ones we intend to pursue.

\section*{Software Availability}

The \texttt{PySymmetry} framework is open-source software distributed under the MIT license. The source code and the issue tracker are publicly available on GitHub. The release that accompanies this paper, v1.0, is archived on Zenodo, and its example notebooks are the ones mapped below to the figures and tables they reproduce.

\begin{itemize}
    \item \textbf{GitHub Repository (Current Version):} \url{https://github.com/ldsufrpe/pysymmetry}
    \item \textbf{Archived release v1.0 (Zenodo DOI):} \url{https://doi.org/10.5281/zenodo.22709708}
    \item \textbf{Documentation:} \url{https://ldsufrpe.github.io/pysymmetry/}
\end{itemize}

The \texttt{Examples} directory of the repository reproduces the results the framework generates, one notebook per artifact. Paths below are relative to \texttt{Examples/}.

\begin{itemize}
    \item \texttt{pde/simulations eigenvalues 2d.ipynb}: Figure~\ref{fig:eigenvalue_comparison}, the correctness validation on the 2D Laplacian.
    \item \texttt{pde/Schrodinger\_benchmark.ipynb}: Table~\ref{tab:performance_schrodinger} and Figures~\ref{fig:performance_schrodinger} and~\ref{fig:speedup_schrodinger}.
    \item \texttt{Chemistry/GF\_problem.ipynb}: the block-diagonal form of Section~\ref{sec:gf} and the listing of \ref{app:GF_example}.
    \item \texttt{Symmetry\_Reduction/D4\_Group.ipynb}: the listing of \ref{app:d4_example}.
\end{itemize}

Two tables are not produced by a notebook. Table~\ref{tab:tool_comparison} is a qualitative comparison compiled from the documentation and the publications of the packages it lists, and Table~\ref{Order_of_groups_and_matrices} reports decompositions obtained in \cite{santos2022decomposition}, where those computations are described in full.
\appendix
\section{Illustrative Workflow: Symmetry Reduction for the \texorpdfstring{$D_4$}{D4} Group}
\label{app:d4_example}

To demonstrate the framework's end-to-end workflow, we consider the dihedral group $D_4$, the symmetry group of a square. Its natural permutation representation acts on a 4-dimensional space corresponding to the vertices of the square. A sample $G$-equivariant matrix $M$ in this space has a structure where the diagonal entries are identical and the off-diagonal entries depend only on whether the vertices are adjacent or diagonally opposite. The complete workflow, from problem definition to solution, is encapsulated in the short script shown in Listing~\ref{lst:d4_workflow}.

Executing this script reveals that the natural representation of $D_4$ decomposes as $\theta = A_1 \oplus B_1 \oplus E$, corresponding to the trivial representation, another one-dimensional representation, and a two-dimensional irreducible representation, respectively. As predicted by the theory, \texttt{PySymmetry} block-diagonalizes the operator $M$ into two $1 \times 1$ blocks and one $2 \times 2$ block, demonstrating the entire algorithmic pipeline on a non-abelian group. For the particular operator used here the two-dimensional block is $9$ times the identity, so the printed result is diagonal throughout.

\begin{lstlisting}[language=Python, caption={A complete workflow in \texttt{PySymmetry} for the symmetry reduction of a $D_4$-equivariant system.}, label={lst:d4_workflow}]
from pysymmetry import FiniteGroup, representation
from sage.all import matrix, QQ

# 1. Define the D4 group using its generators (a rotation and a reflection)
# The group acts on 4 vertices labeled {1, 2, 3, 4}
rotation = "(1,2,3,4)" ; reflection = "(1,3)"
G = FiniteGroup([rotation, reflection], matrix=False)

# 2. Define the natural permutation representation
phi = G.natural_representation()

# 3. Define a G-equivariant matrix M
# 'a' for diagonal, 'b' for adjacent, 'c' for opposite vertices
a, b, c = 10, 2, 1
M = matrix(QQ, [[a, b, c, b],
                [b, a, b, c],
                [c, b, a, b],
                [b, c, b, a]])

# 4. Compute the symmetry-adapted basis via the main function
P = G.base_equivariant_to_blocks(phi)

# 5. Perform the block-diagonalization
M_block_diagonal = P.inverse() * M * P

# Display results
print("Block-Diagonal Matrix P^{-1}MP:\n", M_block_diagonal)

# Block-Diagonal Matrix P^{-1}MP:
#  [15  0  0  0]
# [ 0  7  0  0]
# [ 0  0  9  0]
# [ 0  0  0  9]
\end{lstlisting}

\section{Illustration of the Solution of the GF Problem}
\label{app:GF_example}
\paragraph{Code (\texttt{Sage}/Python using \texttt{pysymmetry}).}
\begin{lstlisting}[language=Python]
from pysymmetry import *
# Defining C2v as a permutation group on its own four elements
# (the regular representation of C2v); the 3x3 matrices below are the
# separate action on the three atoms, so the two degrees differ on purpose.
generators = [[2,1,4,3],[1,2,4,3],[2,1,3,4]]
C2v = FiniteGroup(generators);
gens = [C2v(x) for x in generators ];
# Matrices of permutations of atoms
C2     = matrix(3,3,[[0,1,0],[1,0,0],[0,0,1]]); 
sigmavxz  = matrix(3,3,[[1,0,0],[0,1,0],[0,0,1]]);
sigmavyz  = matrix(3,3,[[0,1,0],[1,0,0],[0,0,1]]);
matrices = [C2, sigmavxz, sigmavyz]
# Defining the representation
rep = representation(gens, matrices);
# Finding the symmetry adapted basis
beta = C2v.base_equivariant_to_blocks(rep); 
beta
# Giving a preview of the block decomposition
print(C2v.quick_block_prevision(rep))
# Defining the matrices F and G as symbolic variables
f11,f12,f33 = var('f11,f12,f33');
g11,g12,g13,g33 = var('g11,g12,g13,g33');
F_H2O = matrix(3,3,[[f11,f12,0],[f12,f11,0],[0,0,f33]]);
G_H2O = matrix(3,3,[[g11,g12,g13],[g12,g11,g13],[g13,g13,g33]]);
# Finding matrices F and G in the new basis
F_new = beta.inverse()*F_H2O*beta;
G_new = beta.inverse()*G_H2O*beta;
FG_new = F_new*G_new;
print(FG_new)

# beta =
# [ 2  0  2]
# [ 2  0 -2]
# [ 0  4  0]
#
# quick_block_prevision: [['degree', 'multiplicity'], [1, 2], [1, 1]]
#
# FG =
# [(f11 + f12)*(g11 + g12)  2*(f11 + f12)*g13                        0]
# [                f33*g13            f33*g33                        0]
# [                      0                  0  (f11 - f12)*(g11 - g12)]
\end{lstlisting}

\bibliographystyle{elsarticle-num}
\bibliography{references}
\end{document}